\documentclass[a4paper, 12pt]{article}
\usepackage[utf8]{inputenc}
\usepackage[top=1.35in, bottom=1.44in, left=1in, right=1in]{geometry}
\usepackage{amsmath, amssymb, amsthm}
\usepackage{mathrsfs}
\usepackage{enumerate}
\usepackage{hyperref}
\usepackage{color} 
\usepackage{float}
\usepackage{xcolor, soul}
\usepackage{tikz}
\usetikzlibrary{arrows.meta, positioning}
\usepackage{parskip}

\newcommand{\nn}{\mathbb{N}}

\newtheorem{definition}{Definition}[section]
\newtheorem{proposition}{Proposition}[section]
\newtheorem{corollary}[proposition]{Corollary}%[section]
\newtheorem{theorem}[proposition]{Theorem}%[section]
\newtheorem{lemma}[proposition]{Lemma}%[section]
\newtheorem{remark}[proposition]{Remark}%[section]

\theoremstyle{definition}
\newtheorem{example}{Example}[section]

\numberwithin{equation}{section}

\newtheorem{theoremA}{Theorem}

\DeclareMathOperator{\diam}{diam}

\title{{\bf Variants of Expansivity, Shadowing and Spectral Decomposition}} 
\author{Pramod Kumar Das$^{a,*}$, Manoj B. Prajapati$^b$\\
\footnotesize $^a$ Department of Mathematics and Statistics, School of Science,\\
\footnotesize Gandhi Institute of Technology and Management (Deemed to be University),\\
\footnotesize Doddaballapur Taluk, Bengaluru, Karnataka 561203, India\\
\footnotesize $^b$ Department of Mathematical Sciences, P. D. Patel Institute of Applied Sciences,\\ 
\footnotesize Charotar University of Science and Technology (CHARUSAT),\\ 
\footnotesize Changa, Anand 388421, India.\\
\footnotesize pramod.math.ju@gmail.com; manoj.prajapati.2519@gmail.com}  
\date{}

\begin{document}
\maketitle

\begin{abstract}
\noindent In this paper, we first distinguish a very recently introduced notion of mean asymptotic expansivity from asymptotic expansivity. Then, we extend the spectral decomposition theorem due to Smale to bi-mean asymptotically expansive homeomorphisms with the shadowing property. We also extend a measurable spectral decomposition theorem due to Shin to shadowable measures for bi-mean asymptotically expansive homeomorphisms. Each of the chain transitive sub-systems obtained by Shin can be further decomposed into chain mixing systems by using the ideas of Richeson and Wiseman. Next, we use a probabilistic argument to recover both Smale-like and Bowen-like spectral decompositions for expansive homeomorphisms with the shadowing property. Moreover, we investigate inter-relations among several variants of the shadowing property including the eventual shadowing property, the thick shadowing property, the ergodic shadowing property, the mean ergodic shadowing property and the $M_{\alpha}$-shadowing property.      
\mbox{}
\bigskip

\mbox{}
\medskip

\end{abstract} 

-------------------------------------

$^*$Corresponding author 

E-mail address: pramod.math.ju@gmail.com (P. K. Das)

\section{Introduction} 

The shadowing property and expansivity are two extensively studied notions in topological dynamics (see \cite{AH}). The interplay between these two notions in the spectral decomposition theorem has been celebrated in topological dynamics for several decades since its inception in the work of Smale in 1967 (see \cite{S1967}) and a subsequent improvement work of Bowen in 1975 (see \cite{B}) for Axiom-A diffeomorphisms on compact manifolds. More specifically, the fact that Anosov diffeomorphisms (a special case of Axiom-A diffeomorphisms) have the shadowing property and expansivity guided Aoki (in 1983) to discard the differentiable structure and prove that the non-wandering set of expansive homeomorphisms with the shadowing property on compact metric spaces admits spectral decomposition (see \cite{AH}). In 1984, Komuro extended it to flows (see \cite{K}). In the first two decades of the twenty first century, many researchers rediscovered, reformulated the spectral decomposition theorem in different settings (see \cite{CDZ},  \cite{CK}, \cite{DLRW}, \cite{JNY}, \cite{LNY}, \cite{LLN}, \cite{O2008}). Very recently, in 2024, Shin obtained a measurable version of the spectral decomposition theorem for expansive shadowable measures on compact metric spaces. One of the purposes of the present paper is to improve and extend this result to shadowable measures for bi-asymptotically expansive and bi-mean asymptotically expansive homeomorphisms on compact metric spaces.    

The notion of asymptotic expansivity was introduced in \cite{LMV} and played an important role in the study of topological spectral decomposition (see \cite[Theorem 9]{LMV}). In \cite{DB}, the authors introduced the notion of mean asymptotic expansivity and investigated its fundamental properties, establishing it as a natural extension of asymptotic expansivity in topological dynamics. 

In this paper, we first construct an example to distinguish mean asymptotic expansivity from asymptotic expansivity. Then, we obtain the following version of Smale's spectral decomposition theorem for bi-mean asymptotically expansive homeomorphisms with the shadowing property.  
\medskip

\begin{theoremA} 
(\textbf{Topological Spectral Decomposition due to Smale})\label{SD}

\noindent Let $X$ be a compact metric space and let $f:X\rightarrow X$ be a bi-mean asymptotically expansive homeomorphism with the shadowing property. Then, there exist finitely many pairwise disjoint closed invariant subsets \(C_1,\ldots,C_l\) of \(\Omega(f)\) such that
\begin{itemize}
\item [(i)] $\Omega(f)=\bigcup_{i=1}^l C_i$,
\item [(ii)] $f|_{C_i}$ is topologically transitive. 
\end{itemize}  
\end{theoremA} 

Using the following result \cite[Theorem 6]{RW}, a new concept called piecewise chain mixing was introduced in \cite{BMR}. Here, we call it category one chain transitivity.  
\medskip

\begin{theorem} 
\label{exclusive} 
For a chain transitive map $f:X\rightarrow X$ of a compact metric space, exactly one of the following holds:   
\begin{itemize}
\item [(1)] there is a period $k\geq 1$ such that $f$ cyclically permutes $k$ closed and open equivalence classes of $X$, and $f^k$ restricted to each equivalence class is chain mixing.
\item [(2)] $f$ factors onto an adding machine map. 
\end{itemize}
\end{theorem} 
\medskip

\begin{definition}
A chain transitive continuous map which satisfies the condition \textit{(1)} of Theorem \ref{exclusive} is called a category one chain transitive map. 
\end{definition} 

We then use this concept to extend \cite[Theorem 1.4]{S} and establish the following measurable version of the spectral decomposition theorem. Consequently, we obtain other versions (see Corollary \ref{Cor1}) one of which is an improvement of \cite[Theorem 1.4]{S}.    
\medskip

\begin{theoremA} (\textbf{Measurable Spectral Decomposition due to Shin and Richeson}) 
\label{SDP} 

Let $\mu$ be a shadowable measure of a bi-mean asymptotically expansive homeomorphism $f$ of an infinite compact metric space $X$. Then,
\begin{itemize}
\item [(1)] There are finitely many disjoint compact invariant chain transitive sets $C_1, \ldots, C_l \subset \Omega(f)$ with $\mu(\Omega(f)) = \sum_{i=1}^{l} \mu(C_i)$. 
\item [(2)] For each $C_j (1\leq j\leq l)$ there is a positive integer $k$ and a finite sequence $B_i$ $(0\leq i\leq k-1)$ of disjoint clopen subsets of $C_j$ such that 

(i) $f(B_i)=B_{{i+1}\pmod{k}}$ and $f^k(B_i)=B_i$, (ii) $C_j=\bigcup_{i=0}^{k-1} B_i$, (iii) $f^k|_{B_i}$ is chain mixing. 
\end{itemize}
\end{theoremA}   

Using Theorem \ref{SDP}, we obtain a probabilistic proof of the following spectral decomposition theorem whose first part is a Smale-like decomposition and the second part is a Bowen-like decomposition (see \cite[Theorem 3.1.11]{AH}). Consequently, we obtain other versions (see Corollary \ref{Cor2}) of this celebrated theorem.    
\medskip

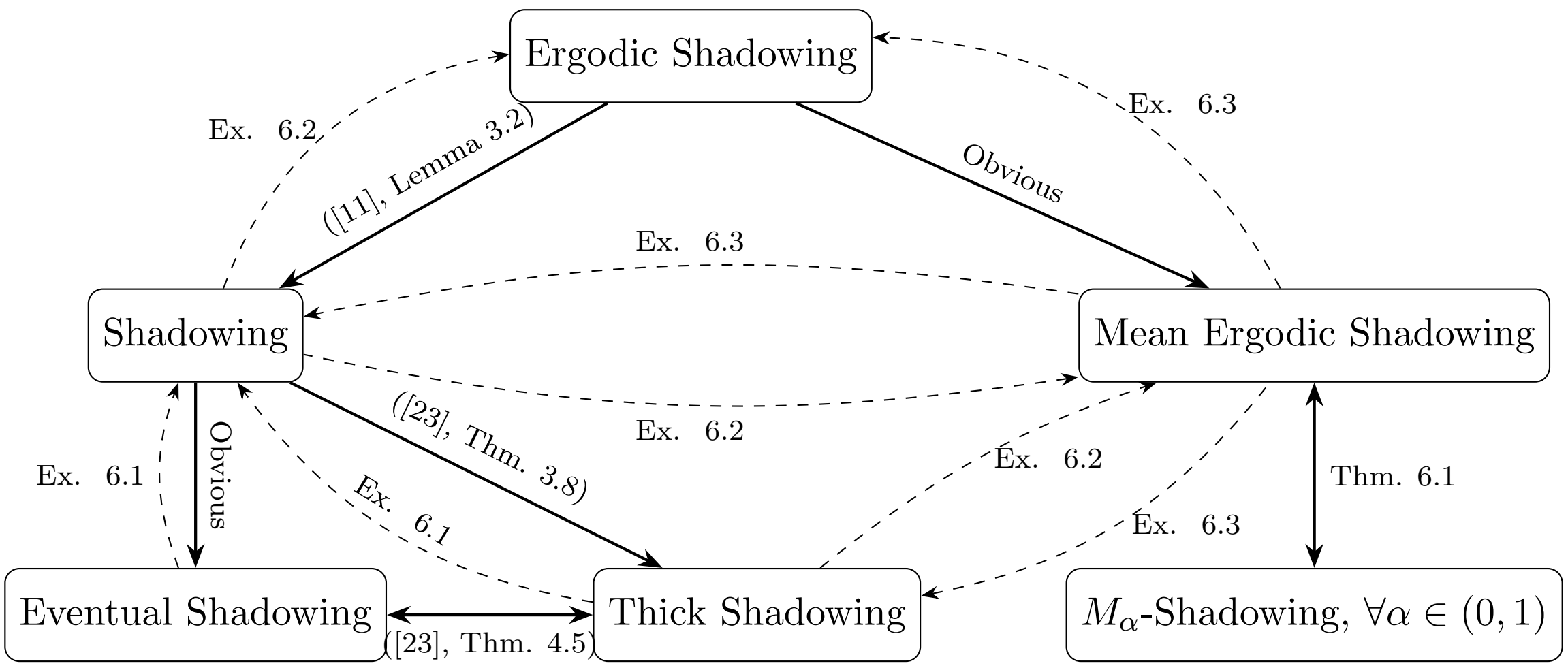
\begin{figure}[H]
\centering
\begin{tikzpicture}[
>=Stealth,
node distance=1.8cm and 2cm,
box/.style={
draw,
rounded corners,
minimum width=2 cm,
minimum height=0.9cm,
align=center,
font=\small
}
]

\node[box] (ES) {Ergodic Shadowing};

\node[box] (S)   [below left=of ES] {Shadowing};

\node[box] (MES) [below right=of ES] {Mean Ergodic Shadowing};

\node[box] (EVS) [below=of S] {Eventual Shadowing};

\node[box] (TS)  [right=of EVS] {Thick Shadowing};

\node[box] (MS) [below=of MES] {$M_{\alpha}$-Shadowing, $\forall \alpha\in (0,1)$}; 

\draw[->, thick]
(ES) -- node[midway, above,sloped] {\scriptsize (\cite{FG}, Lemma 3.2)} (S);

\draw[->, thick]
(ES) -- node[midway, above, sloped] {\scriptsize Obvious} (MES);

\draw[->, thick]
(S) -- node[midway, above, sloped] {\scriptsize Obvious} (EVS);

\draw[->, thick]
(S) -- node[midway, above, sloped] {\scriptsize (\cite{O}, Thm. 3.8)} (TS);

\draw[<->, thick]
(EVS) -- node[below] {\scriptsize (\cite{O}, Thm. 4.5)} (TS);

\draw[<->, thick]
(MES) -- node[right] {\scriptsize Thm. \ref{MES=M_alpha}} (MS);

\draw[->,  dashed]
(S) to[bend left=30]
node[left] {\scriptsize Ex.~ \ref{odometer}}
(ES);

\draw[->, dashed]
(MES) to[bend right=30]
node[right] {\scriptsize Ex.~ \ref{MES but not S}}
(ES);

\draw[->, dashed]
(EVS) to[bend left=20]
node[left] {\scriptsize Ex.~ \ref{ES but not S}}
(S);

\draw[->, dashed]
(TS) to[bend left=22]
node[midway, above, sloped] {\scriptsize Ex.~ \ref{ES but not S}}
(S);

\draw[->, dashed]
(MES) to[bend right=10]
node[above] {\scriptsize Ex.~ \ref{MES but not S}}
(S);

\draw[<-, dashed]
(MES) to[bend left=10]
node[below] {\scriptsize Ex.~ \ref{odometer}}
(S);

\draw[->, dashed]
(TS) to[bend left=10]
node[right] {\scriptsize Ex.~ \ref{odometer}}
(MES);

\draw[<-, dashed]
(TS) to[bend right=20]
node[right] {\scriptsize Ex.~ \ref{MES but not S}}
(MES);
\end{tikzpicture}

\caption{Solid arrows denote implications and dashed arrows denote non-implications} 
\label{shadowing-diagram} 
\end{figure}

\begin{theoremA} (\textbf{Topological Spectral Decomposition due to Shin and Richeson})   
\label{AsymSDP} 

Let $f$ be a bi-mean asymptotically expansive homeomorphism with the shadowing property of an infinite compact metric space $X$. Then,

\begin{itemize}
\item [(1)] There are finitely many disjoint compact invariant transitive sets $C_1, \ldots, C_l \subset \Omega(f)$ with $\Omega(f)=\bigcup_{i=1}^{l} C_i$. 
\item [(2)] For each $C_j (1\leq j\leq l)$ there is a positive integer $k$ and a finite sequence $B_i$ $(0\leq i\leq k-1)$ of disjoint clopen subsets of $C_j$ such that 

(i) $f(B_i)=B_{i+1\pmod{k}}$ and $f^k(B_i)=B_i$, (ii) $C_j=\bigcup_{i=0}^{k-1} B_i$, (iii) $f^k|_{B_i}$ is topologically mixing.    
\end{itemize}   

Motivated by the significant role of the shadowing property in the study of chaos and stability, several variants have been introduced and investigated, including the eventual shadowing property (see \cite{GM},\cite{WY},\cite{XY}), the thick shadowing property (see \cite{DH}, \cite{O}), the ergodic shadowing property (see \cite{FG}), the mean ergodic shadowing property (see \cite{DD}) and the $M_{\alpha}$-shadowing property (see \cite{WOC}). In this paper, we examine and establish the inter-relations among these variants of the shadowing property (see Figure \ref{shadowing-diagram}). More specifically, we prove that the mean ergodic shadowing property is equivalent to the $M_{\alpha}$-shadowing property for all $\alpha\in (0,1)$ (see Theorem \ref{MES=M_alpha}). We provide examples to show that the mean ergodic shadowing property neither implies nor is implied by the thick shadowing property (see Examples \ref{MES but not S} and \ref{odometer}). We provide examples to show that the shadowing property neither implies nor is implied by the mean ergodic shadowing property (see Examples \ref{odometer} and \ref{MES but not S}). It is known that any surjective continuous map with the ergodic shadowing property is chain transitive (see \cite{FG}). We provide an example to show that surjectivity is a necessary condition (see Example \ref{surjectivitynecessary}). It is known that the shadowing property implies the thick shadowing property (see \cite{O}). We construct an example which shows that the compactness of the phase space is a necessary condition (see Example \ref{compactnessnecessary}). We provide an example to show that the thick shadowing property does not imply the shadowing property (see Example \ref{ES but not S}).     
\end{theoremA} 

\section{Definitions and Preliminaries} 

Throughout this paper $(X,d)$ (or simply $X$) denotes a compact metric space unless otherwise stated and $f:X\rightarrow X$ denotes either a homeomorphism or a continuous map. The distance from a point $x\in X$ to a subset $A$ of $X$ is defined by
$\rho(x,A)=\inf \{\, d(x,a)\mid a \in A \,\}$. The $\delta$-neighbourhood of $A$ in $X$ is defined by $B_{\rho}(A,\delta)=\lbrace x\in X\mid \rho(x,A)<\delta\rbrace$. We denote the set of all Borel sets of $X$ by $\mathcal{B}(X)$ and the set of all Borel probability measures of $X$ by $\mathcal{M}(X)$.   
\smallskip

The set of all non-negative integers and the set of all positive integers are denoted by $\mathbb{N}$ and $\mathbb{Z}^+$, respectively. The integer part of a real number $x$ is denoted by $\left\lfloor x\right\rfloor$.   
\smallskip

For $A\subset\mathbb{N}$, the upper density and the lower density of $A$ are respectively defined by 
\noindent $\overline{d}(A)=\limsup_{n\to\infty} \frac{1}{n}|A\cap \lbrace 0,\,1,\ldots,\,n-1\rbrace|$ and $\underline{d}(A)=\liminf_{n\to\infty} \frac{1}{n}|A\cap \lbrace 0,\,1,\ldots,\,n-1\rbrace|$,   
where $|B|$ denotes the number of elements present in the set $B$. If there exists a number $d(A)$ such that
$\overline{d}(A)=\underline{d}(A)=d(A)$, then we say that the set $A$ has density $d(A)$. An infinite set $A\subseteq \mathbb{N}$ with bounded gaps is said to be a syndetic set. An infinite set $A\subseteq \mathbb{N}$ is said to be a thick set if for each $n\in \mathbb{Z}^+$, there exists $i\in\mathbb{N}$ such that $\lbrace i, i+1, i+2, i+3,\cdots,i+n\rbrace\subset A$. A set $A\subseteq \mathbb{N}$ is a thick set if and only if it intersects every syndetic subset of $\mathbb{N}$.   
\smallskip 

Let $f:X\rightarrow X$ be a continuous map. A finite sequence $\lbrace x_i\rbrace_{i=0}^n$ in $X$ is said to be a $\delta$-chain of $f$ if $d(f(x_i),x_{i+1})<\delta$ for $0\leq i<n$. $f$ is said to be chain transitive if for every $\delta>0$ and any pair $x,y\in X$, there exists a $\delta$-chain from $x$ to $y$. $f$ is said to be chain mixing if for every $\delta>0$, there exists a positive integer $N$ such that for any $n\geq N$ and any pair $x,y\in X$, there exists a $\delta$-chain from $x$ to $y$ of length $n$. For $x\in X$ and a pair of non-empty open subsets $U,V$ of $X$, let $N_f(x,U)=\lbrace n\in\mathbb{Z}^+\mid f^n(x)\in U\rbrace$ and $N_f(U,V)=\lbrace n\in\mathbb{Z}^+\mid f^n(U)\cap V\neq \emptyset\rbrace$. A point $x\in X$ is said to be a non-wandering point of $f$ if $N_f(U,U)\neq \emptyset$ for every neighborhood $U$ of $x$. The set of non-wandering points of $f$ is denoted by $\Omega(f)$. $f$ is said to be transitive if for any pair of nonempty open sets $U,V$ in $X$, $N_f(U,V)\neq \emptyset$. $f$ is said to be equicontinuous if for every $\epsilon>0$ there exists $\delta>0$ such that $d(x,y)<\delta$ implies $d(f^n(x),f^n(y))<\epsilon$ for all $n\in\mathbb{N}$.      
\smallskip

Let $f:X\rightarrow X$ be a homeomorphism. A subset $W\subset X$ is called invariant under $f$ if $f(W)=W$. Given $\epsilon > 0$ and $x, y \in X$, we say that $x \sim_{\epsilon} y$ if there is an $\epsilon$-chain from $x$ to $y$ and an $\epsilon$-chain from $y$ to $x$. We say that $x \sim y$ if $x \sim_{\epsilon} y$ for every $\epsilon>0$. The set of chain recurrent points of $f$ is defined as $CR(f)=\lbrace x\in X\mid x \sim x\rbrace$. The equivalence classes of the equivalence relation $\sim$ on $CR(f)$ are called the chain recurrent classes. A point $x\in X$ is called a periodic point with period $N\in\mathbb{Z}^+$ if $f^N(x)=x$. The set of periodic points with period $N$ is denoted by $Fix(f^N)$ and the set of all periodic points is denoted by $Per(f)=\bigcup_{N\in\mathbb{Z}^+} Fix(f^N)$. It is well-known that $CR(f)$ is a closed invariant subset of $X$ and that $Per(f)\subset \Omega(f)\subset CR(f)$.     
\smallskip

Let $f:X\rightarrow X$ be a homeomorphism. Let $p,q\in X$. For $\epsilon>0$, define $N_{\epsilon}^{+}{(p,q)}=\{n\in\mathbb{N}: d(f^n(p),f^n(q))<\epsilon\}$ and ${N}_{\epsilon}^{-}{(p,q)}=\{n\in\mathbb{N}: d(f^{-n}(p),f^{-n}(q))<\epsilon\}$. Then, 

\noindent $\textrm{Asym}_{\epsilon}^+(p)=\Big\lbrace q\in X:  {N}_{\epsilon}^{+}{(p,q)} \text{ is co-finite in } \mathbb{N} \Big\rbrace$,

\noindent $\textrm{Asym}_{\epsilon}^-(p)=\Big\lbrace q\in X:{N}_{\epsilon}^{-}{(p,q)} \text{ is co-finite in } \mathbb{N}\Big\rbrace$, 

\noindent $\textrm{Prox}^+(p)=\Big\lbrace q\in X:$ for every $\epsilon>0,$ $  {N}_{\epsilon}^{+}{(p,q)} \text{ is infinite}\Big\rbrace$, 

\noindent $\textrm{Prox}^-(p)=\Big\lbrace q\in X:$ for every $\epsilon>0,$ ${N}_{\epsilon}^{-}{(p,q)} \text{ is infinite}\Big\rbrace$, 

\noindent $\textrm{ThickProx}^+(p)=\Big\lbrace q\in X:$ for every $\epsilon>0,$ $  {N}_{\epsilon}^{+}{(p,q)} \text{ is thick}\Big\rbrace$,

\noindent $\textrm{ThickProx}^-(p)=\Big\lbrace q\in X:$ for every $\epsilon>0,$ ${N}_{\epsilon}^{-}{(p,q)} \text{ is thick}\Big\rbrace$, 

\noindent $\omega(x)=\Big\lbrace y\in X\mid$ there exists an increasing sequence $\{n_k\}$ in $\nn$ with $\lim\limits_{k\to\infty}f^{n_k}(x)=y\Big\rbrace$,

\noindent $\alpha(x)=\Big\lbrace y\in X\mid$ there exists an increasing sequence $\{n_k\}$ in $\nn$ with $\lim\limits_{k\to\infty}f^{-n_k}(x)=y\Big\rbrace$. 
\medskip

\begin{definition}
Let $f:X\rightarrow X$ be a homeomorphism. A sequence $\lbrace x_i\rbrace_{i\in\mathbb{Z}}$ is said to be a $\delta$-pseudo orbit if $d(f(x_i),x_{i+1})<\delta$ for all $i\in\mathbb{Z}$ and it is said to be $\epsilon$-shadowed by some point $x\in X$ if $d(f^i(x),x_i)<\epsilon$ for all $i\in\mathbb{Z}$. $f$ is said to have the shadowing property if for every $\epsilon>0$ there exists $\delta>0$ such that every $\delta$-pseudo orbit is $\epsilon$-shadowed by some point in $X$.    
\end{definition}
\medskip

\begin{definition}
Let $f:X\rightarrow X$ be a continuous map. A sequence $\lbrace x_i\rbrace_{i\in\mathbb{N}}$ is called a periodic $\delta$-pseudo orbit of $f$ with period $N\in\mathbb{Z}^+$ if $x_{i+N}=x_i$ for all $i\in\mathbb{N}$ and $d(f(x_i),x_{i+1})<\delta$ for all $i\in\mathbb{N}$. $f$ is said to have the strong periodic shadowing property \cite{CDZ} if for every $\epsilon>0$ there exists $\delta>0$ such that for every periodic $\delta$-pseudo orbit $\lbrace x_i\rbrace_{i\in\mathbb{N}}$ of period $N$, there exists $z\in Fix(f^N)$ such that $d(f^i(z),x_i)<\epsilon$ for all $i\in\mathbb{N}$.  
\end{definition}
\medskip

\begin{definition}
Let $f:X\to X$ be a homeomorphism. A point $x\in X$ is
called a \emph{periodic shadowable point} of $f$ if for
every $\varepsilon>0$ there exists $\delta>0$ such that
every periodic $\delta$-pseudo-orbit $(x_i)_{i\in\mathbb Z}$
with $x_0=x$ is $\varepsilon$-shadowed by a periodic point
$p\in\operatorname{Per}(f)$; that is,
$d(f^i(p),x_i)<\varepsilon
\text{ for all }i\in\mathbb Z$. 
We denote the set of all periodic shadowable points of $f$ by
$\mathrm{PSh}(f)
=
\{x\in X : x \text{ is a periodic shadowable point of } f\}$.
\end{definition}
\medskip

\begin{definition}
Let $f:X\to X$ be a continuous map. We say that $f$ has the
\emph{eventual shadowing property} if, for every
$\varepsilon>0$, there exists $\delta>0$ with the following
property: whenever a sequence $(x_n)_{n\ge0}$ satisfies
$d(f(x_n),x_{n+1})<\delta
\text{ for all }n\ge0$, 
there exist a point $z\in X$ and an integer $N\ge0$ such that
$d(f^n(z),x_n)<\varepsilon
\text{ for all }n\ge N$. 
\end{definition}

\section{Asymptotic vs. Mean Asymptotic Expansivity} 

\begin{definition}
Let $f:X\rightarrow X$ be a continuous map. Then, $f$ is said to be asymptotically expansive (see \cite{LMV}) if there exists $\delta>0$ such that $d(f^n(x),f^n(y))\leq \delta$ for all $n\in\mathbb{N}$ implies that $\lim_{n\to \infty} d(f^n(x),f^n(y))=0$. A homeomorphism $f:X\rightarrow X$ is said to be bi-asymptotically expansive if both $f$ and $f^{-1}$ are asymptotically expansive. 
\end{definition}
\medskip

\begin{definition}
Let $f:X\rightarrow X$ be a continuous map. Then, $f$ is said to be mean asymptotically expansive (see \cite{DB}) if there exists $\delta>0$ such that $d(f^n(x),f^n(y))\leq \delta$ for all $n\in\mathbb{N}$ implies that for some $J\subset \mathbb{N}$ with $d(J)=0$, we have $\lim_{n\to \infty} d(f^n(x),f^n(y))=0$ provided $n\notin J$. A homeomorphism $f:X\rightarrow X$ is said to be bi-mean asymptotically expansive if both $f$ and $f^{-1}$ are mean asymptotically expansive. 
\end{definition}

By \cite[Example 4.1]{DB}, asymptotic expansivity implies mean asymptotic expansivity. Therefore, bi-asymptotic expansivity implies bi-mean asymptotic expansivity.  

The following example shows that mean asymptotic expansivity is a weaker notion than asymptotic expansivity. 
\medskip

\begin{example}
Let $\Sigma=\{0,1\}^{\mathbb{N}}$
be endowed with the metric
\[
\rho(u,v)
=
\sum_{i=0}^{\infty}
\frac{|u_i-v_i|}{2^{i+1}}, \text{ where } u=(u_i)_{i\in\mathbb{N}},v=(v_i)_{i\in\mathbb{N}}\in\Sigma.
\] 
Let $\sigma\colon\Sigma\to\Sigma$ denote the one-sided left shift,
defined by $(\sigma u)_i=u_{i+1}$. 
For $z=(z_i)_{i\in\mathbb{N}}\in\Sigma$ given by
\[
z_i
=
\begin{cases}
1, & i=2^r\text{ for some }r\in\mathbb{N},\\
0, & \text{otherwise},
\end{cases}
\]
consider the orbit closure
$Y=\overline{\{\sigma^m z:m\in\mathbb{N}\}}$.  
Then, $Y$ is compact in $\Sigma$ and $\sigma(Y)\subset Y$.
More precisely,
\[
Y
=
\{\sigma^m z:m\in\mathbb{N}\}
\cup\{e_j:j\in\mathbb{N}\}\cup\{\mathbf{0}\},
\]
where $\mathbf{0}=(0,0,\ldots)\in \Sigma$ and $e_j$ is the sequence having exactly one nonzero coordinate, at position $j$. Indeed, for each fixed $j\in\mathbb{N}$, 
$\sigma^{2^r-j}z\rightarrow e_j\text{ as }r\to\infty$. Moreover, by shifting $z$ to the middle of the gaps between consecutive
powers of $2$, one obtains $\mathbf{0}$ as a limit point. Finally, since the gaps between consecutive powers of $2$ tend to infinity, every
other limit point has at most one nonzero coordinate.
\medskip

Now, for $u\in Y$, let $S(u)=\{i\in\mathbb{N}:u_i=1\}$ and
define
\[J_u=\bigcup_{p\in S(u)}
[p-\lfloor\sqrt{p}\rfloor,\,p]\cap\mathbb{N}.
\]
We claim that $\overline{d}(J_u)=0 \text{ for all }u\in Y$. The assertion is immediate when $u=\mathbf{0}$ or $u=e_j$. Suppose
that $u=\sigma^m z$. In this case, $S(u)=\{2^r-m:2^r\geq m\}$.
For \(N\ge2\), choose \(R\in\mathbb N\) such that
$2^R\le 2N+m<2^{R+1}$. If
$[p-\lfloor\sqrt p\rfloor,p]\cap\{0,\ldots,N-1\}\ne\varnothing$, 
then \(p<2N\). Indeed, otherwise \(p\ge2N\ge4\), and hence
$p-\lfloor\sqrt p\rfloor\ge\frac p2\ge N$, 
contradicting the nonempty intersection. Since \(p=2^r-m\), we obtain $2^r=p+m<2N+m<2^{R+1}$, 
so \(r\le R\). Therefore,
\begin{align*} 
&\left|J_u\cap\{0,1,\ldots,N-1\}\right|\leq
\sum_{r=0}^{R}
\left((\sqrt{2})^{r}+1\right)=\sum_{r=0}^{R}
({\sqrt{2}})^{r}+\sum_{r=0}^{R} 1=\frac{(\sqrt{2})^{R+1}-1}{\sqrt{2}-1}+R+1\\&
\leq\frac{\sqrt{2}}{\sqrt{2}-1}(\sqrt{2})^{R}+R+1\leq\frac{\sqrt{2}}{\sqrt{2}-1}\sqrt{2N+m}+\log_2 (2N+m)+1
\end{align*} 
Since \(m\) is fixed, we obtain  
\[
\overline{d}(J_u)
=
\limsup_{N\to\infty}
\frac{
\left|J_u\cap\{0,1,\ldots,N-1\}\right|
}{N}
=0.
\]
Next, we show that 
\[
\lim_{\substack{n\to\infty\\n\notin J_u}}
\rho(\sigma^n u,\mathbf{0})=0.
\]
Suppose that $n\notin J_u$. Then, $u_n=u_{n+1}=\cdots=u_{n+\lfloor\sqrt{n}\rfloor}=0$. 
Otherwise, there would exist $p\in S(u)$ such that
$n\leq p\leq n+\lfloor\sqrt{n}\rfloor$.
Since $p\leq n+\lfloor\sqrt{n}\rfloor$, we have $p-n
\leq\lfloor\sqrt{n}\rfloor
\leq\lfloor\sqrt{p}\rfloor$.
Therefore,
$p-\lfloor\sqrt{p}\rfloor\leq n\leq p$,
which would imply $n\in J_u$, a contradiction.

Therefore, it follows that 
\begin{align*}
\rho(\sigma^n u,\mathbf{0})
&=
\sum_{i=0}^{\infty}
\frac{u_{n+i}}{2^{i+1}}\leq
\sum_{i=\lfloor\sqrt{n}\rfloor+1}^{\infty}
\frac{1}{2^{i+1}}=
2^{-\lfloor\sqrt{n}\rfloor-1} \text{ which gives us } \lim_{\substack{n\to\infty\\n\notin J_u}}
\rho(\sigma^n u,\mathbf{0})=0.  
\end{align*}  
Given $u,v\in Y$, let $J=J_u\cup J_v$. 
Then $\overline{d}(J)=0$, and
\begin{align}\label{uv}
\rho(\sigma^n u,\sigma^n v)
&\leq
\rho(\sigma^n u,\mathbf{0})
+
\rho(\sigma^n v,\mathbf{0})
\rightarrow 0 \text{ as } n\to\infty \text{ with } n\notin J.
\end{align}
Thus, $\sigma|_Y$ is mean asymptotically expansive. On the other hand, since $(\sigma^{2^r}z)_0=1$, we have
$\rho(\sigma^{2^r}z,\mathbf{0})\geq\frac{1}{2} \text{ for all }r\in\mathbb{N}$. In particular, $\rho(\sigma^n z,\mathbf{0})\not\rightarrow 0$.  
\medskip

Now, let us consider   
\[X=\{o\}\cup\bigcup_{k=1}^{\infty}
\left(\{k\}\times(Y\setminus\{\mathbf{0}\})\right).
\]
Define $d$ on $X$ as follows:
\[d((k,u),(\ell,v))=
\begin{cases}
\frac{\rho(u,v)}{2^{k}}, & k=l\\
\frac{\rho(u,\mathbf{0})}{2^{k}} + \frac{\rho(v,\mathbf{0})}{2^{l}} & k\neq l.
\end{cases} \text{ and put } d((k,u),o)=\frac{\rho(u,\mathbf{0})}{2^{k}}.\] 
The triangle inequality follows from the triangle inequality for
$\rho$. For instance,
\[
\frac{\rho(u,v)}{2^{k}}
\leq
\frac{\rho(u,\mathbf{0})}{2^{k}}
+
\frac{\rho(v,\mathbf{0})}{2^{k}}.
\]
Thus, $d$ is a metric on $X$.

We claim that $(X,d)$ is compact. Let $\{x_j\}$ be a sequence in $X$.
If $\{x_j\}$ has infinitely many terms in one fixed copy of $Y$, then
compactness of $Y$ gives a convergent subsequence. Otherwise, after
passing to a subsequence, we may write
$x_j=(k_j,u_j)$ with 
$k_j\rightarrow\infty$.
In this case,
\[
d(x_j,o)
=
\frac{\rho(u_j,\mathbf{0})}{2^{k_j}}
\leq \frac{1}{2^{k_j}}
\rightarrow 0.
\]
Therefore, $x_j\to o$, proving that $X$ is compact.

Now, define $f\colon X\to X$ by $f(o)=o$ 
and $f(k,u)=(k,\sigma u)$, 
where $(k,\mathbf{0})$ is identified with $o$.
For $u,v\in Y$, we have
\[\rho(\sigma u,\sigma v)=
\sum_{i=0}^{\infty}
\frac{|u_{i+1}-v_{i+1}|}{2^{i+1}}\\
=
2\sum_{i=1}^{\infty}
\frac{|u_i-v_i|}{2^{i+1}}\\
\leq 2\rho(u,v).\]
It follows directly from the definition of $d$ that
$d(f(x),f(y))\leq2d(x,y)$ for every $x,y\in X$. 
Thus $f$ is  continuous.

We show that $f$ is mean asymptotically expansive.
Let $x,y\in X$. Suppose first that $x=(k,u)$ and 
$y=(k,v)$ belong to $\{k\}\times Y$. Then
\[
d(f^n(x),f^n(y))
=
\frac{\rho(\sigma^n u,\sigma^n v)}{2^{k}}.
\]
By equation (\ref{uv}), there exists a set $J\subseteq\mathbb{N}$ with
$d(J)=0$ such that
\[
\lim_{\substack{n\to\infty\\n\notin J}}
d(f^n(x),f^n(y))=0.
\]
If $x=(k,u)$ and $y=(\ell,v)$ with $k\neq\ell$, then
\[
d(f^n(x),f^n(y))
=
\frac{\rho(\sigma^n u,\mathbf{0})}{2^{k}}
+
\frac{\rho(\sigma^n v,\mathbf{0})}{2^{l}},
\]
and the same conclusion follows from equation (\ref{uv}). Finally, if $x=(k,u)$ and $y=o$, then $f^n(o)=o$ and
$d(f^n(x),f^n(o))
=2^{-k}\rho(\sigma^n u,\mathbf{0})
\longrightarrow0
\text{ as }n\to\infty,\, n\notin J_u$. Hence $f$ is mean asymptotically expansive.

Finally, we show that $f$ is not asymptotically expansive.
Let $\delta>0$ be arbitrary. Choose $k\geq1$ such that
$\frac{1}{2^{k}}\leq\delta$.
Consider $x=o$ and $y=(k,z)$.
For every $n\in\mathbb{N}$,
\[
d(f^n(x),f^n(y))
=
\frac{\rho(\sigma^n z,\mathbf{0})}{2^{k}}
\leq \frac{1}{2^{k}}
\leq\delta.
\]
Nevertheless, for every $r\in\mathbb{N}$,
\begin{align*}
d(f^{2^r}(x),f^{2^r}(y))
&=\frac{\rho(\sigma^{2^r}z,\mathbf{0})}{2^{k}}
\geq\frac{1}{2(2^{k})}>0.
\end{align*}
It follows that $d(f^n(x),f^n(y))\not\rightarrow 0$. 
Therefore, $f$ is not asymptotically expansive. 
\end{example}

\section{Spectral Decomposition due to Smale} 

The purpose of this section is to extend \cite[Theorem 9] {LMV} to a weaker version of bi-asymptotically expansive homeomorphisms called bi-mean asymptotically expansive homeomorphisms.
\medskip

To prove the following improvement of \cite[Lemma 28]{LMV}, one can use ideas similar to what is used in \cite[Theorem 4.8]{DB}.    
\medskip

\begin{lemma}
\label{SPS}
Let $f:X\rightarrow X$ be a mean asymptotically expansive map. If $f$ has the shadowing property, then it has the strong periodic shadowing property.
\end{lemma} 
\begin{proof}
Let \(c>0\) be a mean asymptotic expansivity constant. Given \(\varepsilon>0\), choose
$0<\eta<\min\{\varepsilon,c/2\}$. 
By the shadowing property, there exists \(\delta>0\) such that every \(\delta\)-pseudo-orbit is \(\eta\)-shadowed by some point of \(X\).

Let \((x_i)_{i\ge0}\) be a periodic \(\delta\)-pseudo-orbit with period \(N\ge1\). Choose \(z\in X\) satisfying
$d(f^i(z),x_i)<\eta$ for all $i\ge0$.
Since \(x_{i+N}=x_i\), we have
\[
\begin{aligned}
d\bigl(f^i(z),f^i(f^N(z))\bigr)
&\le d(f^i(z),x_i)
   +d(x_{i+N},f^{i+N}(z))\\
&<2\eta<c
\end{aligned}
\]
for every \(i\ge0\). Mean asymptotic expansivity therefore gives a set \(J\subseteq\mathbb N\) of density zero such that
\[
\lim_{\substack{i\to\infty\\i\notin J}}
d(f^i(z),f^{i+N}(z))=0.
\]
The set \(N\mathbb N\) has density \(1/N>0\). Consequently, \(N\mathbb N\setminus J\) is infinite, and we can choose an increasing sequence
$n_k\in N\mathbb N\setminus J$. 
By compactness, after passing to a subsequence, we may assume that
$f^{n_k}(z)\rightarrow p \text{ for some }p\in X$. 
Continuity of \(f^N\) yields $d(p,f^N(p))
=\lim_{k\to\infty}
d(f^{n_k}(z),f^{n_k+N}(z))
=0$. 
Thus \(p\in\operatorname{Fix}(f^N)\).

Finally, fix \(i\ge0\). Since each \(n_k\) is a multiple of \(N\), periodicity gives \(x_{n_k+i}=x_i\). Hence
\[
\begin{aligned}
d(f^i(p),x_i)
&=\lim_{k\to\infty}d(f^{n_k+i}(z),x_{n_k+i})\\
&\le\eta<\varepsilon.
\end{aligned}
\]
Therefore \(p\) \(\varepsilon\)-shadows the given pseudo-orbit and satisfies \(f^N(p)=p\). Since \(\delta\) depends only on \(\varepsilon\), this proves the strong periodic shadowing property. 
\end{proof}

Because of Lemma \ref{SPS}, one can follow the same procedure used in \cite[Corollary 29]{LMV} to improve it in the form of the following corollary.   
\medskip

\begin{corollary} \label{CR(f)=Omega(f)}   
Let $f:X\rightarrow X$ be a mean asymptotically expansive map. If $f$ has the shadowing property, then $CR(f)=\Omega(f)=\overline{Per(f)}$.
\end{corollary}  
\medskip

\begin{lemma} 
\label{delta0} 
Let $f:X\rightarrow X$ be a bi-mean asymptotically expansive homeomorphism with the shadowing property. There is $\delta_0>0$ such that $\textrm{ThickProx}^{-}(p)\cap \textrm{ThickProx}^{+}(q)\neq \emptyset$ for all $p,q\in X$ with $d(p,q)<\delta_0$.
\end{lemma} 

\begin{proof}
Let $\varepsilon$ be a bi-mean asymptotically expansive constant of $f$. Further, let $0<\delta_0<\frac{\varepsilon}{2}$ be given for this $\frac{\varepsilon}{2}$ by the shadowing property of $f$. Let $p,q\in X$ with $d(p,q)<\delta_0$. Define the sequence $\{x_n\}_{n\in\mathbb{Z}}$ by 
\[x_i=
\begin{cases}
f^i(p) &\text{ for } i<0\\
f^i(q) &\text{ for } i\geq 0.
\end{cases}
\]
Clearly $\{x_n\}$ is a $\delta_0$-pseudo orbit of $f$. So, there exists $z\in X$ such that $d(f^i(z),x_i)<\frac{\varepsilon}{2}$ for all $i\in\mathbb{Z}$. In other words, $d(f^i(z),f^i(p))<\frac{\varepsilon}{2}$ for all $i<0$ and $d(f^i(z),f^i(q))<\frac{\varepsilon}{2}$ for all $i\geq 0$. Further, $d(z,p)\leq d(z,q)+d(p,q)<\frac{\varepsilon}{2}+\frac{\varepsilon}{2}=\epsilon$. Now, since $\varepsilon$ is a bi-mean asymptotically expansive constant of $f$, there are $J_1,J_2\subset \mathbb{Z}^+$ with $d(J_1)=d(J_2)=0$ such that $\lim_{n\to \infty} d(f^n(z),f^n(q))=0$ provided $n\notin J_1$ and $\lim_{n\to -\infty} d(f^n(z),f^n(p))=0$ provided $n\notin -J_2$. Therefore, $z\in \textrm{ThickProx}^{-}(p)\cap \textrm{ThickProx}^{+}(q)$.  
\end{proof} 

\begin{remark}
\label{remark-delta0}
There is $\delta_0>0$ such that if $p,q\in \Omega(f)$ with $d(p,q)<\delta_0$, then\\ $\textrm{ThickProx}^{-}(p)\cap \textrm{ThickProx}^{+}(q)\cap \Omega(f)\neq \emptyset$. This follows from the proof of the above lemma and the fact that $f|_{\Omega(f)}$ has the shadowing property if $f$ has the same $($see \cite[Theorem 3.1.8]{AH}$)$.  
\end{remark}
\medskip

\begin{lemma}
\label{clopen-shadowing}
Let $f\colon X\to X$ be a homeomorphism of a compact metric space $X$ and let $B$ be a nonempty $f$-invariant clopen subset of $X$. If $f$ has the shadowing property, then
$f|_B\colon B\to B$ also has the shadowing property.  
\end{lemma}

\begin{proof}
The case $B=X$ is immediate. Suppose that $B\neq X$. Since $B$ is
clopen in $X$, the sets $B$ and $X\setminus B$ are disjoint nonempty
compact subsets of $X$. Hence,
\[
r=\operatorname{dist}(B,X\setminus B)=\inf\{d(x,y):x\in B,y\in X\setminus B\}>0. 
\]

Fix $\varepsilon>0$ and put
$\eta=\frac12\min\{\varepsilon,r\}$.
By the shadowing property of $f$, choose $\delta>0$ such that every
$\delta$-pseudo-orbit of $f$ is $\eta$-shadowed. Let $(x_n)_{n\in\mathbb{Z}}$ be a $\delta$-pseudo-orbit of $f|_B$. It is also a $\delta$-pseudo-orbit of $f$, so there exists $z\in X$ such that 
\[
d(f^n(z),x_n)<\eta \text{ for all  } n\in\mathbb{Z}.
\]
Now, $d(z,x_0)<\eta<r$. Since $x_0\in B$ and every point
of $X\setminus B$ is at distance at least $r$ from $x_0$, it follows
that $z\in B$. As $f(B)=B$, we obtain
\[
d\bigl((f|_B)^n(z),x_n\bigr)
=d(f^n(z),x_n)<\eta<\varepsilon
\text{ for all } n\in \mathbb{Z}. 
\]
Thus, $f|_B$ has the shadowing property.   
\end{proof}

\noindent \textbf{Proof of Theorem \ref{SD}}: For $x,y\in X$, we define $x\sim y$ if for every $\epsilon>0$ there exist $\epsilon$-chains from $x$ to $y$ and vice-versa. It is clear that $\sim$ is an equivalence relation on $CR(f)$. Therefore, $CR(f)=\bigcup_{i\in I} C_i$, where $C_i$ ($i\in I$) are disjoint equivalence classes of $\sim$. \textcolor{black}{One can verify that each $C_i$ is compact and $f(C_i)=C_i$}. By Corollary \ref{CR(f)=Omega(f)}, we have $\Omega(f)=\bigcup_{i\in I} C_i$.  
\medskip

Let $\delta_0$ be given by Lemma \ref{delta0}. For $i\in I$, we define
\[B_{\rho}(C_i,\delta_0):=\lbrace p\in \Omega(f): \rho(p,C_i)<\delta_0\rbrace.\] 

We now claim that $B_{\rho}(C_i,\delta_0)\cap Per(f)\subseteq C_i$. 

Take $p \in B_{\rho}(C_i,\delta_0) \cap \text{Per}(f)$. Since $C_i$ is compact, there exists $y \in C_i$ such that $d(p,y) < \delta_0$. Then, the choice of $\delta_0$ and Lemma \ref{delta0} imply that there are
\[
z_1 \in \textrm{ThickProx}^{-}(p) \cap \textrm{ThickProx}^{+}(y) \quad \text{and} \quad z_2 \in \textrm{ThickProx}^{+}(p) \cap \textrm{ThickProx}^{-}(y).
\]
Let $\epsilon>0$. Since $p$ is periodic, there exists $n\in\mathbb{Z}^+$ such that $f^{n}(p)=p$. Since $z_1 \in \mathrm{ThickProx}^{-}(p)$, there exist $k_0\in\mathbb{N}$ and $1 \leq n_1\leq n$ such that
\[
d\big(f^{n_1}(p), f^{-k_0}(z_1)\big) < \epsilon.
\] 
Since $z_1 \in \mathrm{ThickProx}^{+}(y)$ and $C_i$ is an equivalence class under $\sim$, there exists $k_1\in\mathbb{Z}^+$ such that $d\big(f^{k_1}(z_1), f^{k_1}(y)\big)<\epsilon$. Also since $y, f^{k_1}(y)\in C_i$, there exists an $\epsilon$-chain from $f^{k_1}(y)$ to $y$.
Combining these estimates, we obtain that
\[
\{p, f(p), \dots, f^{n_1-1}(p), f^{-k_0}(z_1), f^{-k_0+1}(z_1), \dots, f^{k_1-1}(z_1), f^{k_1}(y),\cdots,y\}
\]
forms an $\epsilon$-chain from $p$ to $y$.  
Now, since $z_2 \in \mathrm{ThickProx}^{-}(y)$ and $z_2 \in \mathrm{ThickProx}^{+}(p)$, there exist $l_1,l_2\in\mathbb{Z}^+$ with $0 \leq n_2 < n$ such that
\[
d\big(f^{-l_1}(y), f^{-l_1}(z_2)\big) < \epsilon
\quad \text{and} \quad
d\big(f^{l_2}(z_2), f^{n_2}(p)\big) < \epsilon.
\]
Since $y, f^{-l_1-1}(y) \in C_i$ and $C_i$ is an equivalence class under $\sim$, there exists an $\epsilon$-chain from $y$ to $f^{-l_1-1}(y)$.
Combining these estimates, we obtain that
\[
\{y,\dots, f^{-l_1-1}(y), f^{-l_1}(z_2), f^{-l_1+1}(z_2), \dots, f^{l_2-1}(z_2), f^{n_2}(p),f^{n_2+1}(p) \dots, p\}
\]
forms an $\epsilon$-chain from $y$ to $p$. Since $\epsilon>0$ was arbitrary,  there exist $\epsilon$-chains from $p$ to $y$ and $y$ to $p$ for all $\epsilon>0$. Thus $p \sim y$, so $p \in C_i$. Therefore, $B_{\rho}(C_i,\delta_0)\cap \mathrm{Per}(f) \subseteq C_i$. 
\medskip

Now we show that $B_{\rho}(C_i,\delta_0) = C_i$ for each $i$. Take any $z \in B_{\rho}(C_i,\delta_0)$. Since by Corollary \ref{CR(f)=Omega(f)}, $\Omega(f) = \overline{\text{Per}(f)}$, there is a sequence $p_k \in \text{Per}(f)$ with $p_k \to z$. So, for large $k$, $p_k \in B_{\rho}(C_i,\delta_0)$ and since $B_{\rho}(C_i,\delta_0)\cap Per(f)\subseteq C_i$, we have $p_k \in C_i$. Since $C_i$ is closed in $\Omega(f)$ and $p_k \to z$, we  get $z\in C_i$ which shows that $B_{\rho}(C_i,\delta_0) \subset C_i$. This means that $B_{\rho}(C_i,\delta_0)=C_i$ and so $C_i$ is open in $\Omega(f)$ for every $i\in I$. Since $\{C_i : i \in I\}$ is an open covering of the compact set $\Omega(f)$, we obtain $\Omega(f) = \bigcup_{i=1}^l C_i$, finite disjoint union of invariant closed subsets $C_i$.    
\medskip

To prove that $f|_{C_i}$ is topologically transitive, first note
that $f|_{\Omega(f)}$ has the shadowing property
(see \cite[Theorem 3.1.8]{AH}). Since $C_i$ is an $f$-invariant
clopen subset of $\Omega(f)$, Lemma \ref{clopen-shadowing} implies that
$f|_{C_i}$ has the shadowing property. Let $U$ and $V$ be nonempty open subsets of $C_i$, and choose
$x\in U$ and $y\in V$. There exist $\varepsilon_1,\varepsilon_2>0$
such that
$B(x,\varepsilon_1)\cap C_i\subseteq U$ and 
$B(y,\varepsilon_2)\cap C_i\subseteq V$. 
Choose $0<\varepsilon<\min\{\varepsilon_1,\varepsilon_2\}$,
and let $\delta>0$ be given by the shadowing property
of $f|_{C_i}$ for this $\varepsilon$.
Since $f|_{C_i}$ is chain transitive, there exists a
$\delta$-chain $x=x_0,x_1,\ldots,x_n=y$ 
contained in $C_i$, with $n\geq 1$.
Extend this chain to a bi-infinite $\delta$-pseudo-orbit
in $C_i$ using the backward orbit of $x$ and the forward
orbit of $y$. By the shadowing property of $f|_{C_i}$,
there exists $z\in C_i$ such that
$d(z,x)<\varepsilon
\quad\text{and}\quad
d(f^n(z),y)<\varepsilon$. 
Since $z,f^n(z)\in C_i$, it follows that $z\in U$ and
$f^n(z)\in V$. Hence $f^n(U)\cap V\neq\varnothing$,
proving that $f|_{C_i}$ is topologically transitive. 

\section{Spectral Decomposition due to Shin and Richeson}    

In this section, we extend and improve a measurable spectral decomposition theorem obtained by Shin (see \cite{S}). We also use a probabilistic argument to recover a topological spectral decomposition theorem.           
\medskip

\begin{definition}
Let $f:X\rightarrow X$ be a homeomorphism. $\mu\in \mathcal{M}(X)$ is said to be shadowable for $f$ (see \cite{S2019}) if for every $\epsilon>0$ there are $\delta>0$ and $B\in \mathcal{B}(X)$ with $\mu(B)=1$ such that every $\delta$-pseudo orbit through $B$ can be $\epsilon$-shadowed by some point $x\in X$.   
\end{definition}          
\medskip

\begin{lemma}[\cite{S}, Lemma 3.5]\label{NW=CR}
If $\mu\in \mathcal{M}(X)$ is a shadowable measure for $f$, then $\mu(\Omega(f))=\mu(CR(f))$.    
\end{lemma}
\medskip

\begin{definition}
Let $f : X \to X$ be a homeomorphism. Then, $\mu\in \mathcal{M}(X)$ is said to have the spectral decomposition property (see \cite{S}) (resp. the category one spectral decomposition property) if there are finitely many disjoint compact invariant sets $C_1, \ldots, C_l \subset \Omega(f)$ with $\mu(\Omega(f)) = \sum_{i=1}^{l} \mu(C_i)$ such that $f|_{C_i}$ is chain transitive (resp. category one chain transitive).     
\end{definition} 
\medskip

\begin{remark}
\label{categoryonedecomposition} 
It is clear from the definition that for a category one chain transitive continuous map $f:X\rightarrow X$ there exists a positive integer $k$ and a finite sequence $B_i$ $(0\leq i\leq k-1)$ of disjoint clopen subsets of $X$ such that  

\noindent (i) $f(B_i)=B_{i+1\pmod{k}}$ and $f^{k}(B_i)=B_i$, (ii) $X=\bigcup_{i=0}^{k-1} B_i$, and (iii) $f^k|_{B_i}$ is chain mixing.   
\end{remark}
\medskip

\begin{definition}
Let $\mu\in \mathcal{M}(X)$ and let $R \subseteq X$ be a  measurable set.
A measurable set $A \subseteq R$ is called $\mu$-open in $R$ if there is $O \subset X$ open such that $O \supset A$ and  $\mu(O \cap R) = \mu(A)$.  
\end{definition}

The following sufficient condition for the spectral decomposition property was established in \cite{S}. 
\medskip

\begin{lemma}[\cite{S}, Lemma 3.3]\label{sufficient} 
Let $f:X \to X$ be a homeomorphism and $\mu\in \mathcal{M}(X)$. If the following two conditions hold, then $\mu$ has the spectral decomposition property. 
\begin{itemize}
\item[(1)] $\mu(\Omega(f)) = \mu(CR(f))$.
\item[(2)] For every chain recurrent component $C$ of $f$, there exists a closed subset $B$ with $\mu(B)=1$ such that $C \cap B$ is $\mu$-open in $CR(f)$. 
\end{itemize} 
\end{lemma} 
\medskip

\begin{lemma}\label{stable-unstable}
Let $C$ and $C'$ be chain recurrent classes of a homeomorphism. If there are $a\in C$ and $a'\in C'$ such that $\textrm{Prox}^+(a)\cap \textrm{Prox}^-(a')\neq\emptyset \text{ and } \textrm{Prox}^-(a)\cap \textrm{Prox}^+(a')\neq\emptyset$, then $C=C'$.
\end{lemma}

\begin{proof}
Let $a\in C$ and $a'\in C'$ such that $\textrm{Prox}^+(a)\cap \textrm{Prox}^-(a')\neq\emptyset$ and $\textrm{Prox}^-(a)\cap \textrm{Prox}^+(a')\neq\emptyset$. We can choose 
$x \in \textrm{Prox}^+(a)\cap \textrm{Prox}^-(a')$ and $x' \in \textrm{Prox}^-(a)\cap \textrm{Prox}^+(a')$. 
\smallskip

Take $b \in \alpha(a)$, $b'\in\alpha(a')$, $d \in \omega(a)$ and $d'\in \omega(a')$. Since $C$ and $C'$ are compact invariant sets, we have $b,d \in C$ and $b',d' \in C'$. 

Fix $\lambda > 0$. Since $f$ is uniformly continuous, there is a $\lambda'>0$ with $\lambda'<\lambda$ such that $d(x,y)<\lambda'$ implies $d(f(x),f(y))<\lambda$. Since $b \in \alpha(a)$, $x' \in \textrm{Prox}^-(a)\cap \textrm{Prox}^+(a')$ and $d'\in \omega(a')$, there exist $n_1,n_2,n_3,n_4\in\mathbb{Z}^{+}$ with $n_1-1>n_2$ and $n_3<n_4$ such that 

\begin{center}
$d(b,f^{-n_1}(a))<\lambda',\,\, d(f^{-n_2}(a),f^{-n_2}(x'))<\lambda,\,\, d(f^{n_3}(x'),f^{n_3}(a'))<\lambda,\,\, d(f^{n_4}(a'),d')<\lambda.$     
\end{center}

So, we can construct a $\lambda$-chain from $b$ to $d'$. Similarly, using $x$, we can construct a $\lambda$-chain from $b'$ to $d$. 
Since $b,d\in C$ and $b',d'\in C'$, we obtain $\lambda$-chain from $d$ to $b$ as well as from $d'$ to $b'$. By concatenating these chains together, we construct a $\lambda$-chain from $b$ to $b'$ as well as from $b'$ to $b$. As $\lambda$ is arbitrary, we have $b \sim b'$. Since
$b\in C$ and $b'\in C'$, we conclude that $C=C'$.
\end{proof} 

\begin{lemma}
\label{chaintoonechain} 
Let \(f\) be a bi-mean asymptotically expansive homeomorphism of an infinite compact metric space \(X\) and let \(\mu\) be a shadowable Borel probability measure for \(f\). If \(\mu\) has the spectral decomposition property, then it has the category one spectral decomposition property. 
\end{lemma}

\begin{proof}
If $\mu(\Omega(f))=0$, take the empty decomposition. Otherwise,
by the spectral decomposition property, there are pairwise disjoint
compact invariant chain transitive sets
$C_1,\ldots,C_l\subseteq\Omega(f)$ such that
\begin{align}
\label{decomposition}  
\mu(\Omega(f))=\sum_{i=1}^{l}\mu(C_i).     
\end{align}
Discard the sets of zero measure. Each remaining $C_i$ lies in
a single chain recurrent class of $f$. Let $K_1,\ldots,K_s$ be
the distinct classes so obtained. They are pairwise disjoint,
compact, invariant and chain transitive, and each has positive
$\mu$-measure. 
	
Let $c>0$ be a common mean asymptotic expansivity constant for
$f$ and $f^{-1}$. By shadowability of $\mu$, there are $\delta_0>0$
and a Borel set $B$ of full measure such that every
$\delta_0$-pseudo-orbit through $B$ is $c/4$-shadowed.
By \cite[Theorem~4.8]{DB}, $\mu$ is periodic shadowable. Hence, by \cite[proof of Lemma~3.5]{DB},
$\mu(\operatorname{PSh}(f))=1$. 
	
Fix $K=K_j$ and choose $a\in K\cap B\cap\operatorname{PSh}(f)$. 
By uniform continuity, choose $0<r<\min\{\delta_0,c/4\}$ such that
\begin{equation}\label{Eq-a}
d(x,y)<r \text{ implies } d(f(x),f(y))<\delta_0.   
\end{equation}

We claim that
\begin{equation}\label{leq-3}
B(a,r)\cap CR(f)\subseteq K.    
\end{equation}

Let \(b\in CR(f)\) satisfy \(d(a,b)<r\). Consider the sequences 
\[
\xi_n=
\begin{cases}
f^n(a),&n\le0,\\
f^n(b),&n>0,
\end{cases}
\qquad
\zeta_n=
\begin{cases}
f^n(b),&n<0,\\
f^n(a),&n\ge0.
\end{cases}
\]
By (\ref{Eq-a}), we have 
\[
d(f(\xi_0),\xi_1)=d(f(a),f(b))<\delta_0.
\]
Also 
\[
d(f(\zeta_{-1}),\zeta_0)=d(b,a)<r<\delta_0.
\]
Therefore, both are \(\delta_0\)-pseudo-orbits through \(a\).
Choose \(u,v\in X\) satisfying
\begin{equation}\label{eq-1}
d(f^n(u),\xi_n)<c/4 \text{ and }
d(f^n(v),\zeta_n)<c/4 \text{ for all } n\in\mathbb{Z}.  
\end{equation}
This implies that 
$d(u,b)\le d(u,a)+d(a,b)<c/4+r<c/2$, 
and similarly \(d(v,b)<c/2\). Together with (\ref{eq-1}), this gives, for every \(n\ge0\),
\[d(f^{-n}(u),f^{-n}(a))<c,\,
d(f^n(u),f^n(b))<c,\,
d(f^n(v),f^n(a))<c,\,
d(f^{-n}(v),f^{-n}(b))<c.\]
    
By bi-mean asymptotic expansivity of $f$, we obtain 
\[
u\in\operatorname{Prox}^{-}(a)\cap\operatorname{Prox}^{+}(b),
\text{ and }
v\in\operatorname{Prox}^{+}(a)\cap\operatorname{Prox}^{-}(b).
\]
Since \(a\) and \(b\) are chain recurrent, Lemma \ref{stable-unstable} shows that their chain recurrent classes coincide. Hence \(b\in K\). That means $B(a,r)\cap CR(f)\subseteq K$.

Fix $x\in K$ and $\varepsilon>0$. Choose
$0<\eta<\min\{r,\varepsilon\}$. Since $a\in\operatorname{PSh}(f)$,
there exists $\delta>0$ such that every periodic $\delta$-pseudo-orbit
through $a$ is $\eta$-shadowed by a periodic point. Concatenate
$\delta$-chains from $a$ to $x$ and from $x$ to $a$, and repeat
the loop to obtain a periodic pseudo-orbit $(y_n)_{n\in\mathbb Z}$
with $y_0=a$ and $y_m=x$. If $x=a$, use a positive-length return
chain and take $m=0$. Let $p$ be a periodic point that $\eta$-shadows this pseudo-orbit.
Then
$d(p,a)<\eta<r$ and 
$d(f^m(p),x)<\eta<\varepsilon$. 
Since $p$ is chain recurrent, \eqref{leq-3} gives
$p\in K$. Thus $f^m(p)$ is a periodic point in $K$ with $d(f^m(p),x)<\varepsilon$. Since $K$ is closed, we obtain
$K=\overline{\operatorname{Per}(f)\cap K}$. Since \(\Omega(f)\) is closed in $X$ and contains every periodic point, $K\subseteq\Omega(f)$. Therefore, each $K_j$ is chain transitive and contains a periodic point. By \cite[Corollary~10]{RW}, $f|_{K_j}$ admits a finite clopen cyclic decomposition on whose pieces the return maps are chain mixing.
    
Finally,
\[
\bigcup_{i=1}^{l}C_i\subseteq\bigcup_{j=1}^{s}K_j
\subseteq\Omega(f).
\]
Then, by (\ref{decomposition}) and the disjointness of the $K_j$ $(1\leq j\leq s)$, we have
\[
\mu(\Omega(f))=\sum_{j=1}^{s}\mu(K_j).
\]
Thus, the sets $K_1,\ldots,K_s$ give the category one spectral
decomposition for $\mu$.    
\end{proof}

\noindent \textbf{Proof of Theorem \ref{SDP}}: Let $\mu$ be a shadowable measure of a bi-mean asymptotically expansive homeomorphism $f : X \to X$ of a compact metric space. Since $\mu(\Omega(f))=\mu(CR(f))$ by Lemma \ref{NW=CR}, it suffices to find $B$ satisfying item (2) of Lemma \ref{sufficient}. 

By the definition of bi-mean asymptotic expansivity, there exists $\epsilon > 0$ such that 
\begin{align*}
\textrm{Asym}^+_{\epsilon}(z) \subset \textrm{Prox}^+(z) 
\quad \text{and} \quad 
\textrm{Asym}^-_{\epsilon}(z) \subset \textrm{Prox}^-(z), 
\quad \forall z \in X.\ldots\ldots\ldots (*)   
\end{align*}

For this $\epsilon$, by the shadowability of $\mu$, we choose $\delta>0$ and $D\in \mathcal{B}(X)$ with $\mu(D)=1$ such that every $\delta$-pseudo orbit through $D$ is $\frac{\epsilon}{2}$-shadowed by some point in $X$. By \cite[Lemma 2.1]{S2019}, every $\delta$-pseudo orbit through $B=\overline{D}$ is $\epsilon$-shadowed by some point in $X$. Now, we prove that this $B=\overline{D}$ satisfies condition \textit{(2)} in Lemma \ref{sufficient}. 
\smallskip

By uniform continuity of $f$, there exists $\delta'$ with $0<\delta'<\delta$ such that $d(x,y)<\delta'$ implies that $d(f(x),f(y))<\delta$. 
\smallskip

For a chain recurrent class $C$, let $B_{\rho}(C,\delta')$ be the $\delta'$-neighborhood of $C$ in $CR(f)$.  If we can prove that $B_{\rho}(C\cap B, \delta')\subset C$, then because of $\mu(B)=1$, we must have $\mu(B_{\rho}(C\cap B,\delta'))\leq \mu(C)=\mu(C \cap B)$. 
This proves that $C \cap B$ is $\mu$-open in $CR(f)$. 
\smallskip

Take $p \in B_{\rho}(C\cap B,\delta')$.  Then there exists $y \in C \cap B$ such that $d(y,p) < \delta'$. 
Define $\lbrace x_i\rbrace_{i\in\mathbb{Z}}$ by 
$x_i=f^i(y)$ if $i\leq 0$ and $x_i=f^i(p)$ if $i> 0$. It follows that $\lbrace x_i\rbrace_{i\in \mathbb{Z}}$ is a $\delta$-pseudo orbit with $x_0 = y \in B$. So, there exists $x \in X$ such that $d(f^i(x), x_i)<\epsilon$ for all $i\in \mathbb{Z}$.   
\smallskip

If $i\leq 0$, then $x \in \textrm{Asym}^-_{\epsilon}(y)$ and if $i>0$, then $x\in \textrm{Asym}^+_{\epsilon}(p)$. So,
$x\in\textrm{Asym}^-_{\epsilon}(y) \cap \textrm{Asym}^+_{\epsilon}(p)$. 
By the choice of $\epsilon$ and the equation $(*)$, we get $x\in \textrm{Prox}^-(y) \cap \textrm{Prox}^+(p)$. This means that $\textrm{Prox}^-(y) \cap \textrm{Prox}^+(p) \neq \emptyset$. 
\smallskip

Repeating the argument but with the sequence $\lbrace y_i\rbrace_{i \in \mathbb{Z}}$ defined by $y_i=f^i(p)$ if $i<0$ and $y_i=f^i(y)$ if $i\geq 0$, we get
$\textrm{Prox}^-(p) \cap \textrm{Prox}^+(y) \neq \emptyset$. 
\smallskip

Therefore, by Lemma \ref{stable-unstable}, we have $p\in C$. 
\smallskip

Thus, $\mu$ has the spectral decomposition property. Then, by Lemma \ref{chaintoonechain}, $\mu$ has the category one spectral decomposition property, that is, there are finitely many disjoint compact invariant sets $C_1, \ldots, C_l \subset \Omega(f)$ with $\mu(\Omega(f)) = \sum_{i=1}^{l} \mu(C_i)$ such that $f|_{C_i}$ is category one chain transitive. Finally, by Remark \ref{categoryonedecomposition},

\noindent (1) There are finitely many disjoint compact invariant chain transitive sets $C_1, \ldots, C_l \subset \Omega(f)$ with $\mu(\Omega(f)) = \sum_{i=1}^{l} \mu(C_i)$.

\noindent (2) For each $C_j (1\leq j\leq l)$ there is a positive integer $k$ and a finite sequence $B_i$ $(0\leq i\leq k-1)$ of disjoint clopen subsets of $C_j$ such that 

\noindent (i) $f(B_i)=B_{i+1\pmod{k}}$ and $f^k(B_i)=B_i$, (ii) $C_j=\bigcup_{i=0}^{k-1} B_i$, (iii) $f^k|_{B_i}$ is chain mixing. 
\medskip

\begin{corollary}
\label{Cor1} 
Every shadowable measure of an (bi-asymptotically) expansive homeomorphism of an infinite compact metric space has the category one spectral decomposition property.   
\end{corollary} 
\medskip

\noindent \textbf{Proof of Theorem \ref{AsymSDP}}: Since the homeomorphism has the shadowing property, every Borel probability measure is shadowable. Since the homeomorphism is bi-mean asymptotically expansive, by Theorem \ref{SDP}, every Borel probability measure has the spectral decomposition property. Then, by \cite[Theorem 1.2]{S}, $f$ has the spectral decomposition property, i.e. there are finitely many non-empty disjoint compact invariant sets $C_1,C_2,\cdots,C_l\subset \Omega(f)$ with $\Omega(f)=\cup_{i=1}^l C_i$ such that $f|_{C_i}$ is chain transitive. We show that $C_j\cap\textrm{Per($f$)}\neq\emptyset$. Choose \(x\in C_j\). Since \(C_j\) is open in \(\Omega(f)\), there exists \(r>0\) such that $B(x,r)\cap\Omega(f)\subseteq C_j$. By Corollary \ref{CR(f)=Omega(f)}, we get $\Omega(f)=\overline{\operatorname{Per}(f)}$. Since \(x\in\overline{\operatorname{Per}(f)}\), $\emptyset\neq B(x,r)\cap\operatorname{Per}(f)\subseteq C_j$. Hence, $C_j$ contains a periodic point. Finally, by \cite[Corollary~10]{RW}, we obtain: 
\noindent (1) There are finitely many disjoint compact invariant chain transitive sets $C_1, \ldots, C_l \subset \Omega(f)$ with $\Omega(f)=\bigcup_{i=1}^{l} C_i$.

\noindent (2) For each $C_j (1\leq j\leq l)$ there is a positive integer $k$ and a finite sequence $B_i$ $(0\leq i\leq k-1)$ of disjoint clopen subsets of $C_j$ such that 

\noindent (i) $f(B_i)=B_{i+1\pmod{k}}$ and $f^k(B_i)=B_i$, (ii) $C_j=\bigcup_{i=0}^{k-1} B_i$, (iii) $f^k|_{B_i}$ is chain mixing. 

Since \(f|_{\Omega(f)}\) has the shadowing property and $C_j$ is a clopen $f$-invariant subset of  \(\Omega(f)\), therefore by Lemma \ref{clopen-shadowing}, \(f|_{C_j}\) has the shadowing property. As $f|_{C_j}$ is chain transitive, it is topologically transitive. Further, by \cite[Remark 2.3.1 and Theorem 2.3.3]{AH}, \(f^k|_{C_j}\) has the shadowing property. Since \(B_i\) is clopen in \(C_j\) and \(f^k(B_i)=B_i\), therefore by Lemma \ref{clopen-shadowing}, \(f^k|_{B_i}\) has the shadowing property. As \(f^k|_{B_i}\) is chain mixing, \cite[Remark 1.5]{Ka24} implies that \(f^k|_{B_i}\) is topologically mixing.  
\medskip

\begin{corollary}
\label{Cor2} 
For every (bi-asymptotically) expansive homeomorphism with the shadowing property on an infinite compact metric space 

\begin{itemize}
\item [(1)] There are finitely many disjoint compact invariant transitive sets $C_1, \ldots, C_l \subset \Omega(f)$ with $\Omega(f)=\bigcup_{i=1}^{l} C_i$. 
\item [(2)] For each $C_j (1\leq j\leq l)$ there is a positive integer $k$ and a finite sequence $B_i$ $(0\leq i\leq k-1)$ of disjoint clopen subsets of $C_j$ such that 

(i) $f(B_i)=B_{i+1\pmod{k}}$ and $f^k(B_i)=B_i$, (ii) $C_j=\bigcup_{i=0}^{k-1} B_i$, (iii) $f^k|_{B_i}$ is topologically mixing.   
\end{itemize}
\end{corollary}

\section{Inter-relations Among Variants of Shadowing}  

Let $f:X\rightarrow X$ be a continuous map. For $x\in X$, $\xi=\lbrace x_i\rbrace_{i\in\mathbb{N}}$ and $\epsilon>0$, define $B(x,\xi,\epsilon)=\lbrace i\in\mathbb{N}\mid d(f^i(x),x_i)<\epsilon\rbrace$. A sequence $\lbrace x_i\rbrace_{i\in\mathbb{N}}$ is said to be a $\delta$-pseudo orbit (resp. $\delta$-ergodic pseudo orbit) of $f$ if $d(f(x_i),x_{i+1})<\delta$ for all $i\in\mathbb{N}$ (resp. $\lbrace i\in\mathbb{N}\mid d(f(x_i),x_{i+1})\geq \delta\rbrace$ has density zero). A sequence $\lbrace x_i\rbrace_{i\in\mathbb{N}}$ in $X$ is said to be $\epsilon$-shadowed (resp. $\epsilon$-ergodic shadowed) (resp. lower $\epsilon$-ergodic shadowed with index $\alpha>0$) (resp. upper $\epsilon$-ergodic shadowed with index $\alpha>0$) (resp. $\epsilon$-thick shadowed) if there exists $x\in X$ such that $d(f^i(x),x_i)<\epsilon$ for all $i\in\mathbb{N}$ (resp. $\lbrace i\in\mathbb{N}\mid d(f^i(x),x_i)\geq \epsilon\rbrace$ has density zero) (resp. {$\underline{d}(B(x,\xi,\epsilon))>\alpha$}) (resp. {$\overline{d}(B(x,\xi,\epsilon))>\alpha$}) (resp. $\lbrace i\in\mathbb{N}\mid d(f^i(x),x_i)<\epsilon\rbrace$ is a thick set).

$(a)$ $f$ is said to have the shadowing property if for every $\epsilon>0$ there exists $\delta>0$ such that every $\delta$-pseudo orbit is $\epsilon$-shadowed. 

$(b)$ $f$ is said to have the ergodic shadowing property (\cite{FG})  if for every $\epsilon>0$ there exists $\delta>0$ such that every $\delta$-ergodic pseudo orbit is $\epsilon$-ergodic shadowed.

$(c)$ $f$ is said to have the mean ergodic shadowing property (\cite{DD})  if for every $\varepsilon\in(0,1)$ there exists $\delta>0$ such that every $\delta$-ergodic pseudo-orbit $\xi=\{x_i\}_{i\in\mathbb N}$ is lower $\varepsilon$-ergodic shadowed with index $1-\varepsilon$ i.e., there exists $x\in X$ such that
$\underline d\bigl(B(x,\xi,\varepsilon)\bigr)
>1-\varepsilon$. 

$(d)$ $f$ is said to have the  $M_{\alpha}$-shadowing property (\cite{DD})  if for every $\epsilon>0$ there exists $\delta>0$ such that every $\delta$-ergodic pseudo orbit is lower $\epsilon$-ergodic shadowed with index $\alpha$. 

$(e)$ $f$ is said to have the  thick shadowing property (\cite{BMR}, \cite{O}) if for every $\epsilon>0$ there exists $\delta>0$ such that every $\delta$-ergodic pseudo orbit is  $\epsilon$-thick shadowed. 

Since the complement of a set with density zero is thick, the ergodic shadowing property implies the thick shadowing property. 
\smallskip 

We first give an example to show that the thick shadowing property does not imply the shadowing property. 
\medskip

\begin{example}\label{ES but not S}
The function $f:[-1,1]\rightarrow[-1,1]$ given by 
\[f(x)=
\begin{cases}
x^3 &\text{ for } -1\leq x\leq 0\\
2x &\text{ for } 0\leq x\leq \frac{1}{2}\\
2(1-x) &\text{ for } \frac{1}{2}\leq x\leq 1
\end{cases}
\]
\end{example}
has the eventual shadowing property (see \cite{GM}). So, by \cite[Theorem 4.5]{O}, $f$ has the thick shadowing property. On the other hand, the authors of \cite{MR} proved that $f$ does not have the shadowing property. 
\medskip

\begin{example}
\label{odometer} 
In \cite[Theorem 3.8]{L}, the author claimed that a chain transitive continuous map with the shadowing property on an infinite compact metric space is sensitive. But an odometer is chain transitive and has the shadowing property. Since the odometer is equicontinuous (see \cite[Remark]{M}, p. 2236), it cannot be sensitive. Therefore, the odometer is an example which shows that \cite[Theorem 3.8]{L} is not true. 
\smallskip

By \cite[Theorem 3.8]{O}, the odometer has the thick shadowing property. Since the ergodic shadowing property implies the mean ergodic shadowing property and the odometer is a minimal map, therefore, \cite[Lemma 3.2]{FG} and \cite[Theorem 3.19]{DD} imply that the odometer cannot have the mean ergodic shadowing property and hence, cannot have the ergodic shadowing property.  
\end{example}

We now characterize the mean ergodic shadowing property in terms of sub-shadowing (see \cite{DH}). 
Although the following equivalence follows from \cite[Theorem 1.2]{DBarXiv}, we prove it here for self-containment.      
\medskip

\begin{theorem}
\label{MES=M_alpha}
A continuous map $f:X\rightarrow X$ has the $M_{\alpha}$-shadowing property for every $\alpha\in (0,1)$ if and only if it has the mean ergodic shadowing property. 
\end{theorem} 

\begin{proof}
Suppose that $f$ has the $M_{\alpha}$-shadowing property, for every $\alpha\in (0,1)$. For $\epsilon_0$ with $0<\epsilon_0<1$ we write $\epsilon_0=1-\alpha_0$ for some $\alpha_0$ with $0<\alpha_0<1$. Let $\delta_0>0$ be given by the $M_{\alpha_0}$-shadowing property of $f$. Let $\xi=\lbrace x_i \rbrace_{i\in\mathbb{N}}$ be a $\delta_0$-ergodic pseudo orbit of $f$. By the $M_{\alpha_0}$-shadowing property of $f$, there exists $x\in X$ such that $\underline{d}(B(x,\xi,\epsilon_0))>\alpha_0$. Therefore, $\overline{d}(B^c(x,\xi,\epsilon_0))=1-\underline{d}(B(x,\xi,\epsilon_0))<1-\alpha_0=\epsilon_0$.  
So, $f$ has the mean ergodic shadowing property.    
\smallskip

Conversely, suppose that $f$ has the mean ergodic shadowing property. For $\epsilon_0$ with $0<\epsilon_0<1$ we write $\epsilon_0=1-\alpha_0$ for some $\alpha_0$ with $0<\alpha_0<1$. Let $\delta_0>0$ be given for this $\epsilon_0$ by the mean ergodic shadowing property of $f$. Let $\xi=\lbrace x_i\rbrace_{i\in\mathbb{N}}$ be a $\delta_0$-ergodic pseudo orbit of $f$. So, there exists $x\in X$ such that $\overline{d}(B^c(x,\xi,\epsilon_0))<\epsilon_0$. Then, we have $\underline{d}(B(x,\xi,\epsilon_0))=1-\overline{d}(B^c(x,\xi,\epsilon_0))>1-\epsilon_0=\alpha_0$. If $\epsilon>\epsilon_0$, then $\underline{d}(B(x,\xi,\epsilon))\geq\underline{d}(B(x,\xi,\epsilon_0))>\alpha_0$. If $\epsilon<\epsilon_0$, by the mean ergodic shadowing property of $f$ there is $\delta>0$ such that every $\delta$-ergodic pseudo orbit $\xi= \lbrace x_i\rbrace_{i\in\mathbb{N}}$ satisfies that $\overline{d}(B^c(x,\xi,\epsilon))<\epsilon$ for some $x\in X$. So, $\underline{d}(B(x,\xi,\epsilon))=1-\overline{d}(B^c(x,\xi,\epsilon))>1-\epsilon>1-\epsilon_0=\alpha_0$. So, $f$ has the $M_{\alpha}$-shadowing property for every $\alpha\in (0,1)$. 
\end{proof} 

\begin{example}
\label{MES but not S}
Let $\mathbb{S}^1$ be the unit circle. Let us consider the surjective continuous map $f:\mathbb{S}^1\rightarrow \mathbb{S}^1$ given by $f(e^{2\pi i\theta})=e^{2\pi i\theta^2}$ (where $\theta\in [0,1)$). In \cite[Example 3.13]{KO}, the authors proved that $f$ has the asymptotic average shadowing property. Then, by \cite[Corollary 3.10]{DD}, $f$ has the mean ergodic shadowing property. Further, by \cite[Theorem 4.5]{XY}, $f$ does not have the eventual shadowing property (or the thick shadowing property) and hence, does not have the shadowing property.   
\end{example}

\begin{center}
\bf \large For surjective continuous maps on compact metric spaces     
\end{center}

\bigskip

The following example shows that the shadowing property and the ergodic shadowing property can be present simultaneously in a system which is not chain transitive. In other words, surjectivity cannot be dropped from the following result.     
\medskip

\begin{lemma}[\cite{FG}, Lemma 3.1]
Any surjective mapping with the ergodic shadowing property is chain transitive. 
\end{lemma}
\medskip

\begin{example} 
\label{surjectivitynecessary} 
Let us consider $X=\bigcup\limits_{n\in\mathbb{Z}^+} X_n\cup \lbrace (0,0)\rbrace$, where $X_n=\lbrace \frac{1}{n}\rbrace\times[0,\frac{1}{2^{n-1}}]$, as a subspace of $\mathbb{R}^2$. Let $f:X\rightarrow X$ be given by

\begin{center}
$f(\frac{1}{k},x)=(\frac{1}{k+1},\frac{1}{2}x)$ for $k\in\mathbb{Z}^+$ and $f(0,0)=(0,0)$. 
\end{center}  

It is easy to verify that $f$ is not chain transitive. We now show that $f$ has the shadowing property and the ergodic shadowing property. 
\smallskip

Indeed, choose $\delta>0$. Since $\diam(X_n)<\delta$ for $n$ with $n>n_0=\left\lfloor \frac{1}{\delta}\right\rfloor$, the uniform continuity of $f$ guarantees that it has the shadowing property.
\smallskip

Let $\epsilon\in (0,\frac{1}{2})$, choose $N\in \mathbb{Z}^{+}$ with $\frac{1}{N+1}<\frac{\epsilon}{\sqrt{2}}\leq \frac{1}{N}$. 

For $\delta=\frac{1}{N}-\frac{1}{N+1}=\frac{1}{N(N+1)}<\frac{\epsilon}{\sqrt{2}}<\epsilon$, let $\lbrace x_i\rbrace_{i\in\mathbb{N}}$ be a $\delta$-ergodic pseudo orbit with \[E=\lbrace i\in\mathbb{N}\mid d(f(x_i),x_{i+1})\geq \delta\rbrace.\] Observe that $\diam(X_i)=\frac{1}{2^{i-1}}<\epsilon$ for all $i\geq N+1$. If $Y=\bigcup_{n\geq N+1} X_n\cup \lbrace (0,0)\rbrace$, then 
\[\diam(Y)=\sqrt{\left(\frac{1}{2^N}\right)^2+\left(\frac{1}{N+1}\right)^2}<\sqrt{\left(\frac{1}{N+1}\right)^2+\left(\frac{1}{N+1}\right)^2}=\frac{\sqrt{2}}{N+1}<\epsilon.\]
Since $E$ has density zero, the set $\bigcup_{i\in E} \lbrace i+1,i+2,i+3,\cdots,i+N\rbrace$ also has density zero. Consequently, $\lbrace i\in\mathbb{N}\mid x_i\in X\setminus Y \rbrace$ has density zero  and hence, any point in $Y$ ergodically $\epsilon$-shadows the $\delta$-ergodic pseudo orbit $\lbrace x_i\rbrace_{i\in\mathbb{N}}$. Therefore, $f$ has the ergodic shadowing property.    
\end{example}

In \cite[Question 3.3]{BMR}, the authors asked whether the shadowing property implies the thick shadowing property. Theorem 3.8 of \cite{O} provides an affirmative answer when the phase space is compact metric. 
\medskip

However, we construct the following example to show that the implication does not hold if the compactness of the phase space is not available.     
\medskip

\begin{example}
\label{compactnessnecessary} 
Let us consider $X=\bigcup_{n\in\mathbb{N}} (\lbrace n\rbrace\times [0,\frac{1}{2^{n}}])\subset\mathbb{R}^2$ as a  subspace of $\mathbb{R}^2$. Let $f:X\rightarrow X$ be given by 
\[f(k,x)=(k+1,\frac{1}{2}x) \text{ for $k\geq 0$}.\]
Choose $\delta>0$. Since $\operatorname{diam}(X_n)<\delta$ for $n$ with 
$n>n_0=\left\lfloor \frac{1}{\delta}\right\rfloor$, the uniform continuity of 
$f$ on the compact set $X_{0}\cup \cdots \cup X_{n_0}$
shows that $f$ has the shadowing property. 
\medskip

Let $\delta>0$ and $D=\lbrace n^2\mid n\in\mathbb{N}\rbrace$. 
Clearly, $D$ has density zero. Let $x_0\in \lbrace 0\rbrace\times [0,1]$. Now, we construct a $\delta$-ergodic pseudo orbit $\lbrace y_j\rbrace_{j\in\mathbb{N}}$ given by 

\begin{center}
$y_0=f(x_0)$, $y_1=f^2(x_0)$ and for each $n\geq 1$, $y_j=f^{(j-n^2)}(x_0)$ for $n^2< j\leq (n+1)^2$.
\end{center}

One cannot $\epsilon$-thick shadow $\lbrace y_j\rbrace_{j\in\mathbb{N}}$ because for any $z\in X$, $f^j(z)$ increases along the $x$-axis as $j$ goes to $\infty$. This shows that $f$ does not have the thick shadowing property.
\end{example}
\medskip

\begin{theorem}
A chain mixing continuous map with the eventual shadowing property (i.e., the thick shadowing property) has the ergodic shadowing property (hence, the mean ergodic shadowing property).  
\end{theorem} 

\begin{proof}
Suppose that $f$ has the eventual shadowing property. Let $\epsilon>0$ be fixed and let $\delta>0$ be given for $\epsilon$ by the eventual shadowing property of $f$. Let \[\eta=\lbrace x_0,x_1,x_2,x_3,...,x_{m_0};x_{m_0+1},x_{m_0+2},...,x_{m_1};x_{m_1+1},x_{m_1+2},...,x_{m_2};x_{m_2+1},...\rbrace\] be a $\delta$-ergodic pseudo orbit of $f$. Let $N\geq 2$ be given by the chain mixing property of $f$. Then, for any $i\in\mathbb{N}$ we consider a $\delta$-chain from $x_{m_i}$ to $x_{m_i+N-1}$ and replace $\lbrace x_{m_i},x_{m_i+1},x_{m_i+2},...,x_{m_i+N-1}\rbrace$ in $\eta$ by this new $\delta$-chain. This process changes $\eta$ to a $\delta$-pseudo orbit. Since $f$ has the eventual shadowing property, there exists $k\in \mathbb{Z}^+$ such that $\eta$ is $\epsilon$-shadowed by an actual orbit after the term $x_{m_k}$. Since $\lbrace m_i\rbrace_{i\in\mathbb{N}}$ has density zero, $\bigcup_{i=k}^{\infty}\lbrace m_i,m_i+1,m_i+2,...,m_i+N-1\rbrace$ also has density zero.  Thus $f$ has the ergodic shadowing property.
\end{proof}

\begingroup
\setlength{\parskip}{2pt}
\setlength{\itemsep}{0pt}
\setlength{\parsep}{0pt}

\endgroup
\end{document}